\documentclass[english]{article}
\usepackage[T1]{fontenc}
\usepackage[latin9]{inputenc}
\usepackage{fancyhdr}
\usepackage{float}
\usepackage{textcomp}
\usepackage{amsmath}
\usepackage{amssymb}
\usepackage{graphicx}
\usepackage{esint}

\makeatletter

\newcommand{\lyxmathsym}[1]{\ifmmode\begingroup\def\b@ld{bold}
  \text{\ifx\math@version\b@ld\bfseries\fi#1}\endgroup\else#1\fi}

\floatstyle{ruled}
\newfloat{algorithm}{tbp}{loa}
\providecommand{\algorithmname}{Algorithm}
\floatname{algorithm}{\protect\algorithmname}

\usepackage[T1]{fontenc}
\usepackage{pgfplots}
\usepackage{times}
\date{}

\usepackage{lastpage}
\makeatother

\usepackage{babel}
\begin{document}

\title{Estimation of errors in iterative solutions of a non-symmetric linear
system}

\author{Aashish Vishwakarma$^{*}$, Murugesan Venkatapathi$^{+}$\\
Supercomputer~Education~and~Research~Centre\\
 Indian Institute of Science, Bangalore, India\\
 aashish.vishwakarma@gmail.com$^{*}$, murugesh@serc.iisc.ernet.in$^{+}$ }
\maketitle
\begin{abstract}
Estimation of actual errors from the residue in iterative solutions
is necessary for efficient solution of large problems when their condition
number is much larger than one. Such estimators for conjugate gradient
algorithms used to solve symmetric positive definite linear systems
exist. This work presents error estimation for iterative solutions
of general indefinite linear systems to provide accurate stopping
and restarting criteria. In many realistic applications no properties
of the matrices are known a priori; thus requiring such a general
algorithm.

Our method for approximating the required quadratic form $r^{T}A^{-1}r$
(square of the $A-$ norm of the error vector) when solving nonsymmetric
linear systems with Bi-Conjugate Gradient (BiCG) algorithm, needs
only $O(1)$ time (per BiCG iteration). We also extend this estimate
to approximate $l_{2}-$norm of error vector using the relations of
Hestenes and Stiefel \cite{hestenes1952methods}. Using the heuristics
of numerical results we observe that the developed algorithm (BiCGQL)
is at least $\varkappa\times10^{-1}$ times more accurate than residue
vector based stopping criteria (where $\varkappa$ is the condition
number of the system).

\end{abstract}

\section{Introduction}

\label{intro}

When using any iterative algorithm for solving a linear system {[}$Ax=b${]},
one of the most important questions is when to stop the iterations.
One would like to stop the iterations when the norm of the error (where
$x_{k}$ are the iterates) 
\begin{equation}
\epsilon_{k}=x-x_{k}\label{eq:error}
\end{equation}

is small enough. However, the actual error is unknown and most iterative
algorithms rely on residual vector as stopping criteria like $||r_{k}||_{2}\approx||\epsilon||_{2}||b||_{2}$
where $r_{k}=b-Ax_{k}$ is the residual vector. Such stopping criteria
can work when the system is well-conditioned. However, these types
of stopping criteria can be misleading depending on the condition
number of $A$ or the choice of the initial approximation. It can
stop the iterations too early when the norm of the error is still
much larger than tolerance, or too late in which case too many floating
point operations have been done for obtaining the required accuracy. 

Also, in the case when the condition number of the matrix is too large,
the residual vector might not show proper converging behavior. In
fact the residual vector might show oscillating behavior while the
actual error might still be (however slowly) converging (and vice-versa).
In such cases, the residual vector cannot be used as a good stopping
or restarting criteria. The norm of relative residue can be as large
as $\kappa$ times or as small as $1/\kappa$ times the norm of the
relative error. 

Even though most iterative algorithms are used with preconditioner,
in most realistic situations, it is not guaranteed whether preconditioner
will actually reduce the condition number of the matrix (this can
mostly be seen in case of large matrices). This created motivation
for ways to compute estimates of some norms of the error for iterative
solvers. Such estimators (e.g. CGQL) are already available for CG
algorithm (\cite{meurant2006lanczos}). For solving nonsymmetric linear
systems using FOM and GMRES methods formulas for estimation of errors
have been suggested \cite{meurant2011estimates} recently. Our objective
is to derive efficient estimators for solving general nonsymmetric
linear systems. 

In our paper, we briefly recall the error estimates existing for Hermitian
Positive Definite (HPD) problems (section -2). We show the equivalence
conditions between (non symmetric) Lanczos co-efficients and BiCG
iterates (sections 3.2-3.3). We develop efficient $O(1)$ estimations
for $A-$norm and $l_{2}-$norm of the error vector for general indefinite
problems (sections 3.4-3.5) using local bi-orthogonality conditions.
The estimation formulas we derive depend only upon BiCG iterates and
add no extra cost to the BiCG algorithm. We test this method (BiCGQL)
for BiCG computations and compare them with the residual based stopping
criteria and existing bounds for non-symmetric problems suggested
by Golub and Meurant in their book (\cite{golub2009matrices}, p.210).
We show that our estimators result in large improvements of the stopping
criteria as the condition number of the problems increase.

\section{Related Work}

\begin{algorithm}[h]
\caption{Lanczos Algorithm\label{alg:Lanczos-Algorithm}}

\textbf{Input} $A,\: v$

$\beta_{0}=0,\: v_{0}=0$

$v_{1}=v/||v||$

\textbf{for $k=1....$convergence}

$\qquad$$\begin{aligned}w & =Av_{k}-\beta_{k-1}v_{k-1}\\
\alpha_{k} & =v_{k}^{T}w\\
w & =w-\alpha_{k}v_{k}\\
\beta_{k} & =||w||\\
v_{k+1} & =w/\beta_{k}
\end{aligned}
$

\textbf{end for}
\end{algorithm}

One of the most commonly used methods for solving linear systems with
real symmetric positive definite (SPD) matrix is Conjugate Gradient
(CG) algorithm. It can be derived from several different perspectives,
(i) an orthogonalization problem, (ii) minimization problem and (iii)
Lanczos algorithm. In their paper, Golub and Meurant \cite{golub1997matrices}
have suggested computing bounds for A-norm of the error in the Conjugate
Gradient (CG) method. A typical norm of error for CG is the A-norm
(also called the energy norm) which is minimized at each iteration.
It is defined as 
\begin{equation}
||\epsilon_{k}||_{A}^{2}\equiv(\epsilon_{k}^{T}A\epsilon_{k})=r^{T}A^{-1}AA^{-1}r=r^{T}A^{-1}r
\end{equation}

It is sometimes also interesting to compute $l_{2}-$norm, for which
$||\epsilon||_{2}^{2}=r^{T}A^{-2}r.$ In order to obtain bound on
$||\varepsilon||_{A}$ we must obtain bound on $r^{T}A^{-1}r$. Of
course we do not want to compute $A^{-1}$. So our problem is similar
to obtaining computable bound for quadratic forms $u^{T}A^{-1}u.$

In \cite{meurant2006lanczos}, Meurant has showed how one can obtain
approximation for $A-$norm of error in CG iterations. 

When $A=A^{T}$,
\begin{align}
||\varepsilon^{2}||_{A} & =r^{T}A^{-1}r\nonumber \\
 & =r^{T}Q\Lambda^{-1}Q^{T}r\nonumber \\
 & =q^{T}\Lambda^{-1}q\nonumber \\
 & =\sum_{i=1}^{n}\lambda_{i}^{-1}q_{i}^{2}\label{eq:eq3}\\
 & =\int_{{{\lambda}_{\min}}}^{{{\lambda}_{\max}}}{{{\lambda}^{-}}^{1}d\alpha(\lambda)}\label{eq:eq4}\\
 & =\int_{{a}}^{{b}}{f(\lambda)d\alpha(\lambda)}\mbox{ (In general)}\label{eq:eq5}
\end{align}

Equation \ref{eq:eq4} is Riemann\textendash{}Stieltjes integral of
equation \ref{eq:eq3}. Here $\alpha$ is a piecewise constant and
defined as

\begin{align}
\alpha(\lambda) & =0 & \mbox{ if }\lambda\leq\lambda_{min}\\
 & =\sum_{j=1}^{i}q_{j}^{2} & \mbox{ if }\lambda_{i}\leq\lambda<\lambda_{i+1}\\
 & =\sum_{j=1}^{n}q_{j}^{2} & \mbox{\ensuremath{\mbox{ if }\lambda\geq\lambda_{max}}}
\end{align}

This allows us to use Gauss, Gauss-Radau, and Gauss-Lobatto formulas
for a function $f$ given by (from equation \ref{eq:eq5})

\begin{equation}
\int_{a}^{b}f(\lambda)d\alpha(\lambda)=\sum_{i=1}^{N}w_{i}f(t_{i})+\sum_{j=1}^{M}v_{j}f(z_{j})+R[f]\label{eq:integral}
\end{equation}

where the weights $w_{i}$,$v_{j}$and nodes $t_{i}$ are unknowns
and nodes $z_{j}$ are given. $R[f]$ can be given by

\begin{equation}
R[f]=\frac{f(\eta)^{2N+M}}{(2N+M)!}\int_{a}^{b}\prod_{j=1}^{M}(\lambda-z_{j})\left(\prod_{i=1}^{N}(\lambda-t_{i})\right)^{2}d\alpha(\lambda)
\end{equation}

where, $a<\eta<b$

\vspace{5 mm}

When $M=0$, the approximation of the integral (equation \ref{eq:integral})
is called the Gauss formula, when $M=1,\; z_{1}=\lambda_{min}\;\mbox{or}\; z_{1}=\lambda_{max}$
it is called Gauss-Radau and when $M=2,\; z_{1}=\lambda_{min}\mbox{ and }\; z_{2}=\lambda_{max}$
it is called Gauss-Lobatto. The nodes $t$ and $z$ can be obtained
by a polynomial decomposition of the integral in terms of $p_{i}(\lambda)$.
Moreover, a set of orthogonal polynomials provides a 3-term recursion
relationship for easy calculations. This means the recurrence coefficients
can be represented in a matrix of symmetric tri-diagonal form; the
crucial observation being that these can be trivially extracted from
the CG iterates, resulting in negligible addition of computation over
the iterative solution. In more generality, the CG algorithm can be
described as a minimization of the polynomial relation

\begin{equation}
||x-x_{k}||_{A}=\underset{p_{k}}{\mathrm{min}}||p_{k}(A)(x-x_{0})||_{A}
\end{equation}

Given $\int_{a}^{b}p_{i}(\lambda)p_{j}d\alpha(\lambda)=0$ when $i\neq j$
and $1$ when $i=j$ and $\gamma_{i}p_{i}(\lambda)=(\lambda-\omega_{i})p_{i-1}(\lambda)+\gamma_{i-1}p_{i-2}(\lambda),\mbox{ when }i=1....N$,
normalized such that

\[
\int d\alpha=1,\: p_{0}(\lambda)=1,\: p_{-1}(\lambda)=0
\]

\begin{equation}
\Rightarrow\lambda P_{N-1}(\lambda)=T_{N}P_{N-1}(\lambda)+\gamma_{N}p_{N}(\lambda)\varepsilon_{N}
\end{equation}

where $\varepsilon_{N}^{T}=[0,0,0...1],P_{N-1}(\lambda)^{T}=[p_{1}(\lambda),...p_{N-1}(\lambda)]$
and $T_{N}$ is the Jacobi matrix obtained by Lanczos Algorithm as
discussed later. These techniques are used for providing lower and
upper bounds for quadratic forms $u^{T}f(A)u$ where $f$ is a smooth
function, $A$ is an SPD matrix and $u$ is a given vector. Paper
by Golub, Gene H., and Zden\v{e}k Strako� \cite{golub1994estimates}
talks about how to obtain error estimations in quadratic formulas.
The algorithm GQL (Gauss Quadrature and Lanczos) is based on the Lanczos
algorithm and on computing functions of Jacobi matrices. These techniques
are adapted to the CG algorithm to compute lower and upper bounds
on the A-norm of the error. The idea is to use CG instead of the Lanczos
algorithm, to compute explicitly the entries of the corresponding
Jacobi matrices from the CG coefficients, and then to use the same
formulas as in GQL. The formulas are summarized in the CGQL algorithm
(QL standing for Quadrature and Lanczos) (algorithm \ref{alg:CGQL}).
The CGQL algorithm uses the tridiagonal Jacobi matrix obtained by
translating the coefficients computed in CG into the Lanczos coefficients.
This paper focuses on establishing the relationships between (non-symmetric)
Lanczos co-efficients and BiCG iterates in order to obtain expressions
for approximation of $A-$norm of the error vector and further extending
the approach for obtaining approximation of $l_{2}-$norm of the error
vector.

\subsection*{The Lanczos, CG and CGQL Algorithms}

Given a starting vector $v$ and an SPD matrix $A$, the Lanczos algorithm
(algorithm \ref{alg:Lanczos-Algorithm}) computes an orthonormal basis
$v_{1},...,v_{k+1}$ of the Krylov subspace $\kappa_{k+1}(A,v)$,
which is defined as...
\begin{equation}
\kappa_{k+1}(A,v)=span\{v,Av,...A^{k}v\}
\end{equation}

In Algorithm \ref{alg:Lanczos-Algorithm} we have used the modified
Gram-Schmidt form of the algorithm. The basis vectors $v_{j}$ satisfy
the matrix relation

\begin{equation}
AV_{k}=V_{k}T_{k}+\eta_{k+1}v_{k+1}\varepsilon_{k}^{T}
\end{equation}

Here, $\varepsilon_{k}$ is the $k^{th}$ canonical vector, where
$V_{k}=[v_{1}...v_{k}]$ and $T_{k}$ is the $k\times k$ symmetric
tridiagonal matrix of the recurrence coefficients computed in algorithm
\ref{alg:Lanczos-Algorithm}:
\begin{equation}
T_{k}=\left[\begin{array}{cccc}
\alpha_{1} & \eta_{1}\\
\eta_{1} & \ddots\\
 &  & \ddots & \eta_{k-1}\\
 &  & \eta_{k-1} & \alpha_{k}
\end{array}\right]
\end{equation}

The coefficients $\beta_{j}$ being positive, $T_{k}$ is a Jacobi
matrix. The Lanczos algorithm works for any symmetric matrix, but
if $A$ is positive definite, then $T_{k}$ is positive definite as
well.

When solving a system of linear algebraic equations $Ax=b$ with symmetric
and positive definite matrix $A$, the CG method (algorithm \ref{alg:Conjugate-Gradient-Algorithm})
can be used. CG (which may be derived from the Lanczos algorithm)
computes iterates $x_{k}$ that are optimal since the A-norm of the
error defined in (1) is minimized over $x_{0}+\kappa_{k}(A,r_{0})$,
\begin{equation}
||x-x_{k}||_{A}=\underset{y\in x_{0}+\kappa_{k}(A,r_{0})}{\mathrm{min}}||x-y||_{A}
\end{equation}

\begin{algorithm}[H]
\textbf{input} $A,\: b,\: x_{0}$

$r_{0}=b-Ax_{0}$

$p_{0}=r_{0}$

\textbf{for} $k=1...n$ \textbf{do}

$\qquad$$\begin{aligned}\gamma_{k-1} & =\frac{r_{k-1}^{T}r_{k-1}}{p_{k-1}^{T}Ap_{k-1}}\\
x_{k} & =x_{k-1}+\gamma_{k-1}p_{k-1}\\
r & =r_{k-1}-\gamma_{k-1}Ap_{k-1}\\
\beta_{k} & =\frac{r_{k}^{T}r_{k}}{r_{k-1}^{T}r_{k-1}}\\
p_{k} & =r_{k}+\delta_{k}p_{k-1}
\end{aligned}
$

\textbf{end for}

\caption{Conjugate Gradient Algorithm\label{alg:Conjugate-Gradient-Algorithm}}
\end{algorithm}

It is well-known that the recurrence coefficients computed in both
algorithms (Lanczos and CG) are connected via
\begin{equation}
\eta_{k}=\frac{\sqrt{\beta_{k}}}{\gamma_{k-1}},\:\alpha_{k}=\frac{1}{\gamma_{k-1}}+\frac{\beta_{k-1}}{\gamma_{k-2}},\:\delta_{0}=0,\:\gamma_{-1}=1
\end{equation}

Noticing that the error $\epsilon_{k}$  and the residual $r_{k}$
are related through $A\epsilon_{k}=r_{k}$ , we have

\begin{equation}
||\epsilon||_{A}^{2}=\epsilon_{k}^{T}A\epsilon_{k}=r_{k}^{T}A^{-1}r_{k}
\end{equation}

The above formula has been used for reconstructing the A-norm of the
error. For the sake of simplicity, only Gauss rule has been considered.
Let $s_{k}$ be the estimate of $||\varepsilon^{k}||_{A}$. Let $d$
be a positive integer, then the idea is to use the following formula
at CG iteration $k$, 

In their paper \cite{golub1997matrices} , Golub and Meurant give
following expression pertaining to $A-$norm of the error. 

\begin{equation}
||\epsilon||_{A}^{2}=||r_{0}||^{2}((T_{n}^{-1})_{1,1}-(T_{k}^{-1})_{1,1})
\end{equation}

Further, for sufficiently large $k$, and $d=1$, they denote the
estimator as

\begin{equation}
s_{k-1}=||r^{0}||_{2}^{2}\frac{\eta_{k-1}^{2}c_{k-1}}{\delta_{k-1}(\alpha_{k}\delta_{k-1}-\eta^{2})}>0
\end{equation}

Using rules of Gauss, Gauss-Radalu, and Gauss Lobatto error bound,
the CGQL algorithm can be established as algorithm \ref{alg:CGQL}.
It should be noted that, in their more recent work, Gerard Muerant\cite{meurant2005estimates}
has derived formula relating the $l_{2}-$norm of the error in CG
algorithm. 

\begin{algorithm}[H]
\textbf{input} $A,\: b,\: x_{0},\:\lambda_{m},\:\lambda_{M}$ 

$r_{0}=b\lyxmathsym{\textminus}Ax_{0},\:\mbox{ }p_{0}=r_{0}$

$\eta_{0}=0,\:\gamma_{\lyxmathsym{\textminus}1}=1,\: c_{1}=1,\:\beta_{0}=0,\:\delta_{0}=1,\:\bar{\alpha}(\mu)^{1}=\lambda_{m},\:\underline{\alpha}(\eta)^{1}=\lambda_{M}$

\textbf{for} $k=1....$until convergence \textbf{do}

$\qquad$CG-iteration (k)

$\qquad$$\begin{aligned}\alpha_{k} & =\frac{1}{\gamma_{k-1}}+\frac{\beta_{k-1}}{\gamma_{k-2}}, & \eta_{k}^{2} & =\frac{\beta_{k}}{\gamma_{k-1}^{2}}\\
\delta_{k} & =\alpha_{k}-\frac{\beta_{k-1}^{2}}{\delta_{k-1}}, & g_{k} & =||r_{0}||\frac{c_{k}^{2}}{\delta_{k}}\\
\overline{\delta_{k}} & =\alpha_{k}-\overline{\alpha_{k}},\;\overline{\alpha_{k+1}}=\lambda_{m}+\frac{\beta^{2}}{\overline{\delta_{k}}},\qquad & \overline{f_{k}} & =||r_{0}||^{2}\frac{\eta_{k}^{2}c_{k}^{2}}{\delta_{k}(\overline{\alpha_{k+1}}\delta_{k}-\eta_{k}^{2})}\\
\underline{\delta_{k}} & =\alpha_{k}-\underline{\alpha_{k}},\;\underline{\alpha_{k+1}}=\lambda_{M}+\frac{\beta^{2}}{\overline{\delta_{k}}}, & \underline{f_{k}} & =||r_{0}||^{2}\frac{\eta_{k}^{2}c_{k}^{2}}{\delta_{k}(\underline{\alpha_{k+1}}\delta_{k}-\eta_{k}^{2})}\\
\breve{\alpha}_{k+1} & =\frac{\overline{\delta_{k}}\underline{\delta_{k}}}{\overline{\delta_{k}}-\underline{\delta_{k}}}\left(\frac{\lambda_{m}}{\overline{\delta_{k}}}-\frac{\lambda_{M}}{\underline{\delta_{k}}}\right), & \breve{\eta}_{k} & =\frac{\overline{\delta_{k}}\underline{\delta_{k}}}{\overline{\delta_{k}}-\underline{\delta_{k}}}(\lambda_{M}-\lambda_{m})\\
\overline{f_{k}} & =||r_{0}||^{2}\frac{[\breve{\eta}_{k}]^{2}c_{k}^{2}}{\delta_{k}(\breve{\alpha}_{k+1}\delta_{k}-[\breve{\eta}_{k}]^{2})}\\
c_{k+1}^{2} & =\frac{\eta_{k}^{2}c_{k}^{2}}{\delta_{k}^{2}}
\end{aligned}
$

\textbf{end for}

\caption{\label{alg:CGQL}CGQL (Conjugate Gradients and Quadrature via Lanczos
coefficients)}
\end{algorithm}

\section{Methodology}

\subsection{The Problem}

One of the most useful algorithm for iterative solution of non-symmetric
linear systems in context of Lanczos and CG algorithms is Bi-Conjugate
Gradient (Bi-CG) algorithm. 

A-norm of error in BiCG can be written as, 
\begin{equation}
||\epsilon_{A}||^{2}=e^{T}Ae=r^{T}A^{-1}r\label{eq:anorm}
\end{equation}

Here, $r$ is residual vector pertaining to the BiCG method. When
$A$ is positive definite, the right side of the above equation is
always positive, it is also called as energy norm in physics related
problems. In case of indefinite matrices, the absolute value of the
above equation is considered. Moreover, the $l_{2}-$norm of the error
can be written as, 

\begin{equation}
||\epsilon||^{2}=e^{T}e=r^{T}A^{-2}r\label{eq:l2norm}
\end{equation}

We are interested in approximating \ref{eq:anorm} and \ref{eq:l2norm}.
In their paper, Starkov and Tichy \cite{strakovs2011efficient} develop
a method of $O(\sim n)$ to approximate a bilinear form $(c^{T}Ab)$
based on BiCG method. Our goal is to approximate the quantity $r_{k}^{T}A^{-1}r_{k}$
(A-norm of the error) for every iteration of a BiCG algorithm. In
following sections, we will derive the approximation for $A-$norm
and $l_{2}-$norm of the error for every BiCG iteration. The BiCG
method is shown as algorithm \ref{alg:BiCG}.

\begin{algorithm}[H]
\textbf{input:} $A,\: A^{T},\: x_{0},\: b$

$r_{0}=b-Ax_{0}$,$\:$ $\tilde{r_{0}}=p_{0}=\tilde{p_{0}}=r;$

\textbf{for} $k=1....$

$\qquad\begin{aligned}\alpha_{k} & =\frac{\tilde{r}_{k}^{T}r_{k}}{p_{k}^{T}Ap_{k}}\\
x_{k+1} & =x_{k}+\alpha_{k}p_{k},\qquad & \tilde{x}_{k+1} & =\tilde{x}_{k}+\alpha_{k}\tilde{p}_{k}\\
r_{k+1} & =r_{k}-\alpha_{k}Ap_{k}, & \tilde{r}_{k+1} & =\tilde{r}_{k}-\alpha_{k}A^{T}\tilde{p}_{k}\\
\beta_{k+1} & =\frac{\tilde{r}_{k+1}^{T}r_{k+1}}{\tilde{r}_{k}^{T}rk}\\
p_{k+1} & =r_{k+1}+\beta_{k+1}p_{k}, & \tilde{p}_{k+1} & =\tilde{r}_{k+1}+\beta_{k+1}\tilde{r}_{k}
\end{aligned}
$

\textbf{end}

\caption{BiCG Algorithm\label{alg:BiCG}}
\end{algorithm}

\subsection{Non-Symmetric Lanczos algorithm}

Let A be a non-singular matrix of order $n$. We introduce the Lanczos
algorithm as a means of computing an orthogonal basis of a Krylov
subspace. Let $v_{1}$ and $\tilde{v}_{1}$ be given vectors (such
that $||v_{1}||=1$ and $(v_{1},\tilde{v}_{1})=1)$. 

For $k=1,2,...$

\begin{equation}
\begin{aligned}z_{k} & =Av_{k}-\omega_{k}v_{k}-\eta_{k-1}v_{k-1}\\
w_{k} & =A^{T}\tilde{v}_{k}-\omega_{k}\tilde{v}_{k}-\tilde{\eta}_{k-1}\tilde{v}_{k-1}
\end{aligned}
\end{equation}

The coefficient $\omega_{k}$ being computed as $\omega_{k}=(\tilde{v}_{k},Av_{k}).$
The other coefficients $\eta_{k}$ and $\tilde{\eta}_{k}$ are chosen
(provided $(z_{k},w_{k})=0$) such that $\eta_{k}\tilde{\eta}_{k}=(z_{k},w_{k})$
and the new vectors at step $k+1$ are given by

\begin{equation}
\begin{array}{cc}
v_{k+1} & =\frac{z_{k}}{\tilde{\eta}_{k}}\\
\tilde{v}_{k+1} & =\frac{w_{k}}{\eta_{k}}
\end{array}
\end{equation}

These relations can be written in matrix form, let

\begin{equation}
T_{k}=\left(\begin{array}{ccccc}
\omega_{1} & \eta_{1}\\
\tilde{\eta}_{1} & \omega_{2} & \eta_{2}\\
 & \ddots & \ddots & \ddots\\
 &  & \tilde{\eta}_{k-1} & \omega_{k-1} & \eta_{k-1}\\
 &  &  & \tilde{\eta}_{k} & \omega_{k}
\end{array}\right)
\end{equation}

and 

\begin{equation}
\begin{array}{cc}
V_{k} & =[v_{1}...v_{k}]\\
\tilde{V_{k}} & =[\tilde{v}_{1}...\tilde{v}_{k}]
\end{array}
\end{equation}

then

\begin{equation}
\begin{array}{cc}
AV_{k} & =V_{k}T_{k}+\tilde{\eta}_{k}v_{k+1}(\varepsilon_{k})^{T}\\
A^{T}\tilde{V}_{k} & =\tilde{V}_{k}T_{k}^{T}+\eta_{k}\tilde{v}_{k+1}(\varepsilon_{k})^{T}
\end{array}
\end{equation}

In order to approximate $A^{-1}$, which is restricted onto $\kappa_{n}(A,r_{0})$,
following holds for non-symmetric Lanczos algorithm

\begin{equation}
A^{-1}=V_{n}T_{n}^{-1}\tilde{V}_{n}^{T}\label{eq:2}
\end{equation}

Considering the starting vectors $v_{1}=r_{0}/||r_{0}||$ and $w_{1}=||r_{0}||/r_{0}$,
we get

\begin{align}
r_{0}^{T}A^{-1}r_{0} & =\frac{r_{0}^{T}r_{0}}{||r_{0}||}w_{1}V_{n}T_{n}^{-1}\tilde{V}_{n}v_{1}||r_{0}||\nonumber \\
 & =(r_{0}^{T}r_{0})\varepsilon_{1}^{T}T_{n}^{-1}\varepsilon_{1}\label{eq:3}\\
 & =||r_{0}||^{2}(T_{n}^{-1})_{(1,1)}\label{eq:eqa30}
\end{align}

where $\varepsilon_{1}$is first canonical vector. In the next section
we are going to establish relationship between BiCG iterates and Lanczos
co-efficients.

\subsection{BiCG , Gauss Quadrature and Lanczos (BiCGQL)}

For a square matrix $A$, having the distribution function $w(\lambda)$
and interval $(a,b)$ such that $a<\lambda_{1}<\lambda_{2}...<\lambda_{n}<b$,
for any continuous function, one can define Riemann-Stieltjes integral
as 
\begin{equation}
\int_{a}^{b}f(\lambda)dw(\lambda)\label{eq:EQ31}
\end{equation}

where $w(\lambda)$ is a stepwise constant function.

\[
w(\lambda)=\begin{cases}
0 & \mbox{for \ensuremath{\lambda<\lambda_{1}}}\\
\sum\limits _{j=1}^{i}w_{j} & \mbox{for \ensuremath{\lambda_{i}\leq\lambda<\lambda_{i+1}\mbox{, \ensuremath{1\leq i\leq n-1}}}}\\
\sum\limits _{j=1}^{n}w_{j} & \mbox{for \ensuremath{\lambda_{n}>\lambda}}
\end{cases}
\]

Integral \ref{eq:EQ31} is a finite sum,

\[
\int_{a}^{b}f(\lambda)dw(\lambda)=\sum\limits _{i=1}^{n}w_{i}f(\lambda_{i})=v_{1}^{T}f(A)v_{1}
\]

We are interested in the quadratic formula, $r_{k}^{T}A^{-1}r_{k}$.
It can be written using Riemann-Stieltjes integral for function $f(\lambda)=1/\lambda$.
In $n^{th}$ step of non-symmetric Lanczos algorithm we get the full
orthonormal basis of $\kappa_{n}(A,v_{1})$ and we have

\begin{align}
AV_{n} & =V_{n}T_{n}\\
\Rightarrow A^{-1}V_{n} & =V_{n}T_{n}^{-1}
\end{align}

and

\begin{align}
\int_{a}^{b}f(\lambda)dw(\lambda) & =\sum\limits _{i=1}^{n}w_{i}f(\lambda_{i})=v_{1}^{T}f(A)v_{1}\nonumber \\
 & =v_{1}^{T}A^{-1}\varepsilon_{1}=v_{1}^{T}V_{n}(T_{n}^{-1})\varepsilon_{1}\nonumber \\
 & =\varepsilon_{1}^{T}(T_{n}^{-1})\varepsilon_{1}=(T_{n}^{-1})_{1,1}\label{eq:equa34}
\end{align}

From above equation and equation \ref{eq:eqa30}, it can be said that
BiCG can implicitly compute weights and nodes of Gauss Quadrature
rule applied to Riemann-Stieltjes integral as 

\begin{equation}
\int_{a}^{b}f(\lambda)dw(\lambda)=(T_{n}^{-1})_{1,1}=\frac{||x-x_{0}||_{A}^{2}}{||r_{0}||^{2}}\label{eq:equa35}
\end{equation}

As mentioned earlier, using Gauss rule on the interval $[a,b]$ and
a function $f$ (such that its Riemann-Stieltjes integral and all
moments exist), the above function can be approximated as

\begin{equation}
\int_{a}^{b}f(\lambda)dw(\lambda)=\sum\limits _{i=1}^{k}w_{i}f(v_{i})+R_{k}^{G}
\end{equation}

In Lanczos terms, it can be expressed as

\begin{equation}
(T_{n}^{-1})_{1,1}=(T_{k}^{-1})_{1,1}+R_{k}^{G}\label{eq:eq11}
\end{equation}

The remainder is nothing but scaled $A-norm$ of the error.

\begin{equation}
R_{k}^{G}=\frac{r_{k}^{T}A^{-1}r_{k}}{||r_{0}||^{2}}\label{eq:equa 38}
\end{equation}

Using BiCG iterates from algorithm, the relation between $r_{0}$
and $r_{k}$ can be written as

\begin{equation}
r_{0}^{T}A^{-1}r_{0}=\sum\limits _{j=0}^{k}\alpha_{j}||r_{j}||^{2}+r_{k}^{T}A^{-1}r_{k}\label{eq:equa 39}
\end{equation}

for which, the Gauss Quadrature approximation is (using \ref{eq:eqa30},
\ref{eq:eq11}, \ref{eq:equa 38} and \ref{eq:equa 39})

\begin{equation}
(T_{k}^{-1})_{1,1}=\frac{1}{||r_{0}||^{2}}\left(r_{0}^{T}A^{-1}r_{0}-r_{k}^{T}A^{-1}r_{k}\right)=\frac{1}{||r_{0}||^{2}}\sum\limits _{j=0}^{k-1}\alpha_{j}||r_{j}||^{2}\label{eq:eq15}
\end{equation}

\subsection{$O(1)$ expression for approximating $A-$norm of the error}

Let us again consider Gauss Quadrature rule at step $k$. 

\begin{equation}
(T_{n}^{-1})_{1,1}=(T_{k}^{-1})_{1,1}+\frac{r_{k}^{T}A^{-1}r_{k}}{||r_{0}||^{2}}\label{eq:eq16}
\end{equation}

Here, we want to approximate $r_{k}^{T}A^{-1}r_{k}$. Of course at
iteration $k<n$, $(T_{n}^{-1})_{1,1},$ is not known. Re-writing
the above equation at step $k+1$, 

\begin{equation}
(T_{n}^{-1})_{1,1}=(T_{k+1}^{-1})_{1,1}+\frac{r_{k+1}^{T}A^{-1}r_{k+1}}{||r_{0}||^{2}}\label{eq:eq17}
\end{equation}

Subtracting \ref{eq:eq17} from \ref{eq:eq16}, we get

\begin{align}
r_{k}^{T}A^{-1}r_{k}-r_{k+1}^{T}A^{-1}r_{k+1} & =||r_{0}||^{2}[(T_{k+1}^{-1})_{1,1}-(T_{k}^{-1})_{1,1}]\nonumber \\
 & =[\alpha_{k}||r_{k}||^{2}]\mbox{ (from equation \ref{eq:eq15})}\label{eq:eq18}
\end{align}

\ref{eq:eq18} gives insights for approximating $r_{k}^{T}A^{-1}r_{k}$.
Alternatively, the same expression can be derived using the expressions
for $r_{k+1}$ and $p_{k}$, in BiCG algorithm as following \\
for $k=0,1,2,3\mbox{...}n-1$

\begin{align}
r_{k}^{T}A^{-1}r_{k}-r_{k+1}^{T}A^{-1}r_{k+1} & =(r_{k+1}+\alpha_{k}A^{T}p_{k})^{T}A^{-1}(r_{k+1}+\alpha_{k}Ap_{k})-r_{k+1}^{T}A^{-1}r_{k+1}\nonumber \\
 & =(r_{k+1}^{T}+\alpha_{k}p_{k}^{T}A)(A^{-1}r_{k+1}+\alpha_{k}p_{k})-r_{k+1}^{T}A^{-1}r_{k+1}\nonumber \\
 & =r_{k+1}^{T}A^{-1}r_{k+1}+\alpha_{k}p_{k}^{T}r_{k+1}+\alpha_{k}r_{k+1}^{T}p_{k}\nonumber \\
 & +\alpha_{k}^{2}p_{k}^{T}Ap_{k}-r_{k+1}^{T}A^{-1}r_{k+1}\nonumber \\
 & \Rightarrow r_{k}^{T}A^{-1}r_{k}-r_{k+1}^{T}A^{-1}r_{k+1}=\alpha_{k}r_{k}^{T}r_{k}\label{eq:32}
\end{align}

$(\because r_{k+1}^{T}p_{k}=p_{k}^{T}r_{k+1}=0)$

From equations \ref{eq:32}, 

\[
r_{k}^{T}A^{-1}r_{k}-r_{k+1}^{T}A^{-1}r_{k+1}=\alpha_{k}||r_{k}||^{2}
\]

\begin{equation}
\implies||\epsilon_{k}||_{A}^{2}-||\epsilon_{k+1}||_{A}^{2}=\alpha_{k}||r_{k}||^{2}\label{eq:anormrecur}
\end{equation}

(as $r^{T}A^{-1}r=||\epsilon||_{A}^{2}$($A-$norm of the error))

For a finite natural number $d(\geq0)$ the above expression can be
approximated as 

\begin{equation}
||\epsilon_{k-d}||_{A}^{2}\approx\sum\limits _{j=k-d}^{k}\alpha_{j}||r_{j}||^{2}\label{eq:anrmapprox}
\end{equation}

Here, $d$ signifies the delay in approximation. It should be noted
that when $A$ is positive-definite, the above expression is always
positive and thus provides a lower bound for the square of $A-$norm
of the error. When $A$ is indefinite, the above expression can be
negative (upper bound) or positive (lower bound) depending up on residual
vector. Also, the BiCG method might show irregular convergence, in
such cases higher values of $d$ can result can less accurate approximations.
Hence, lower values of $d$ is recommended. Method like BiCGSTAB can
repair the irregular convergence behavior of BiCG, as a result smoother
convergence is obtained, and hence higher values of $d$ can be used
for better approximations.

\subsection{$O(1)$ expression for approximation of $l_{2}$-norm of the error}

Hestenes and Stiefel \cite{hestenes1952methods} proved the following
result relating the $l_{2}-$norm and $A-$norm of the error.

\begin{equation}
||\epsilon_{k}||_{A}^{2}+||\epsilon_{k+1}||_{A}^{2}=[||\epsilon_{k}||^{2}-||\epsilon_{k+1}||^{2}][\mu(p_{k})]\label{eq:41}
\end{equation}

\begin{equation}
\mbox{where, }\mu(p_{k})=\frac{(p_{k}^{T}Ap_{k})}{||p_{k}||^{2}}
\end{equation}

Summing equations \ref{eq:anormrecur} and \ref{eq:41}, 

\begin{align}
2||\epsilon_{k}||_{A}^{2} & =[||\epsilon_{k}||^{2}-||\epsilon_{k+1}||^{2}][\mu(p_{k})]+\alpha_{k}||r_{k}||^{2}\nonumber \\
\implies||\epsilon_{k}||^{2}-||\epsilon_{k+1}||^{2} & =\frac{2||\epsilon_{k}||_{A}^{2}}{[\mu(p_{k})]+\alpha_{k}||r_{k}||^{2}}\label{eq:43}\\
\implies||\epsilon_{k}||^{2}-||\epsilon_{k+1}||^{2} & =\phi_{k},\mbox{ }\left(\mbox{where, \ensuremath{\phi_{k}}=\ensuremath{\frac{2||\epsilon_{k}||_{A}^{2}}{[\mu(k)]+\alpha_{k}||r_{k}||^{2}}}}\right)
\end{align}

For a finite natural number $d$$(\geq0)$ , above expression can
be used to approximate $l_{2}-$norm of the error as following

\begin{equation}
||\epsilon_{k-d}||^{2}\approx\sum\limits _{j=k-d}^{k}\phi_{j}\label{eq:l2norm approx}
\end{equation}

Here, $d$ signifies the delay in approximation. It should be noted
that if $d_{1}$delay is introduced in the estimation of $A-$norm
of error in equation \ref{eq:anrmapprox}, and $d_{2}$delay is introduced
in estimation of $l_{2}-$norm, the total delay becomes $d_{1}+d_{2}$
and the equation \ref{eq:l2norm approx} becomes

\[
||\epsilon_{k-d_{1}-d_{2}}||^{2}\approx\sum\limits _{j=k-d_{1}-d_{2}}^{k-d_{1}}\phi_{j}
\]

These estimator are incorporated with BiCG in algorithm \ref{alg:BiCG-Bound}
(we call it BiCGQL-BiCG Quadrature Lanczos).

\begin{algorithm}[h]
\textbf{input:} $A,A^{T},x_{0},b,d_{1},d_{2}$

$r_{0}=b-Ax_{0}$, $\tilde{r_{0}}=p_{0}=\tilde{p_{0}}=r;$

\textbf{for} $k=0,1....$

$\qquad$$\begin{aligned}\alpha_{k} & =\frac{\tilde{r}_{k}^{T}r_{k}}{p_{k}^{T}Ap_{k}}\\
\mu(p_{k}) & =\frac{(p_{k}^{T}Ap_{k})}{||p_{k}||^{2}}
\end{aligned}
$

$\qquad$$\begin{aligned}x_{k+1} & =x_{k}+\alpha_{k}p_{k}, & \tilde{x}_{k+1} & =\tilde{x}_{k}+\alpha_{k}\tilde{p}_{k}\\
r_{k+1} & =r_{k}-\alpha_{k}Ap_{k},\qquad & \tilde{r}_{k+1} & =\tilde{r}_{k}-\alpha_{k}A^{T}\tilde{p}_{k}
\end{aligned}
$

\textbf{$\qquad$if} $k\geq d_{1}+d_{2}$

$\qquad\qquad\begin{aligned}g_{k-d_{1}} & =\sum\limits _{j=k-d_{1}}^{k}\alpha_{j}||r_{j}||^{2}\mbox{ (\ensuremath{A-}norm estimation)}\\
\phi_{k-d_{1}} & =\frac{2||\epsilon_{k-d_{1}}||_{A}^{2}}{[\mu(p_{k-d_{1}})]+\alpha_{k-d_{1}}||r_{k-d_{1}}||^{2}}\\
f_{k-d_{1}-d_{2}} & =\sum\limits _{j=k-d_{1}-d_{2}}^{k-d_{1}}\phi_{j}\mbox{ (\ensuremath{l_{2}-}norm estimation)}
\end{aligned}
$

$\qquad$\textbf{end if}

$\qquad$$\beta_{k+1}=\frac{\tilde{r}_{k+1}^{T}r_{k+1}}{\tilde{r}_{k}^{T}rk}$

$\qquad$$p_{k+1}=r_{k+1}+\beta_{k+1}p_{k},\quad\quad\tilde{p}_{k+1}=\tilde{r}_{k+1}+\beta_{k+1}\tilde{r}_{k}$

\textbf{end}

\caption{BiCGQL Algorithm\label{alg:BiCG-Bound}}
\end{algorithm}

\subsection{BiCG Convergence}

Few theoretical results are known about the convergence of BiCG. For
HPD systems the method delivers the same results as CG, but at twice
the cost per iteration. For nonsymmetric matrices it has been shown
that in phases of the process where there is significant reduction
of the norm of the residual, the method is more or less comparable
to full GMRES (in terms of numbers of iterations) (Freund and Nachtigal
\cite{freund1991qmr}). In practice this is often confirmed, but it
is also observed that the convergence behavior may be quite irregular,
and the method may even break down. The breakdown situation due to
the possible event can be circumvented by so-called look-ahead strategies
(Parlett, Taylor and Liu \cite{parlett1985look}). The other breakdown
situation, occurs when the decomposition fails, and can be repaired
by using another decomposition (such as QMR developed by Freund and
Nachtigal \cite{freund1991qmr}\cite{freund1994implementation}).
Sometimes, breakdown or near-breakdown situations can be satisfactorily
avoided by a restart at the iteration step immediately before the
(near) breakdown step. BiCGSTAB is an improvement over BiCG algorithm
which leads to a considerably smoother convergence behavior. It should
be noted that relations \ref{eq:anrmapprox} and \ref{eq:l2norm approx}
hold valid for BiCGSTAB algorithms.

\section{Numerical Results and Validations }

\subsection{Tests and Results}

\begin{figure}[h]
\begin{centering}
\includegraphics[scale=0.5]{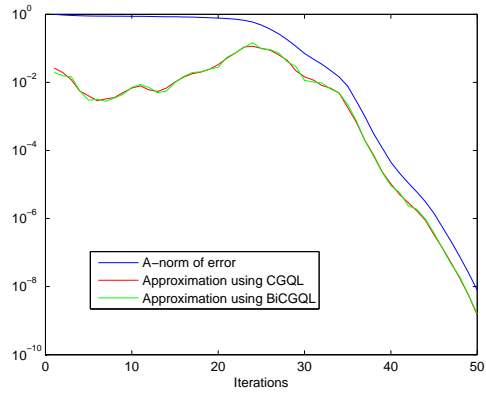}
\par\end{centering}

\caption{Comparing BiCGQL A-norm estimator with Gauss Approximation of CGQL
for an HPD matrix (condition number of the matrix is $10^{4},$ $d_{1}=d_{2}=0$)
\label{fig:HPD-case:-Comparing}}
\end{figure}

In figure \ref{fig:HPD-case:-Comparing}, the developed approximator
is compared with the $A-$norm of the error vector as well as Gauss
Approximation (discussed in CGQL algorithm). Here, $A$ is an square
HPD matrix. Figure shows that our $A-$norm estimator is almost as
good as CGQL Gauss Rule for HPD matrices. 

\begin{figure}[h]
\begin{centering}
\includegraphics[scale=0.5]{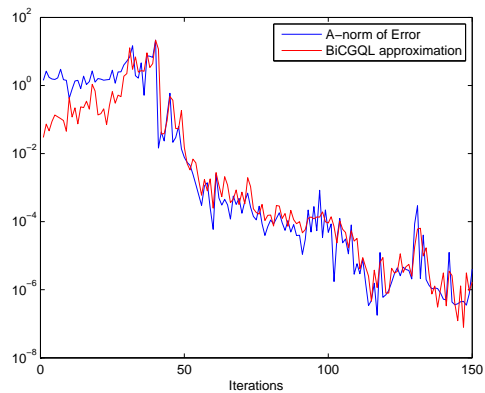}
\par\end{centering}

\caption{BiCGQL estimator in case of a non-hermitian (indefinite) matrix (absolute
values are considered for quadratic term $r^{T}A^{-1}r$ and its approximation,
condition number of matrix is $10^{4}$, $d_{1}=d_{2}=0$)\label{fig:3}}
\end{figure}
\begin{figure}[h]
\begin{centering}
\includegraphics[scale=0.5]{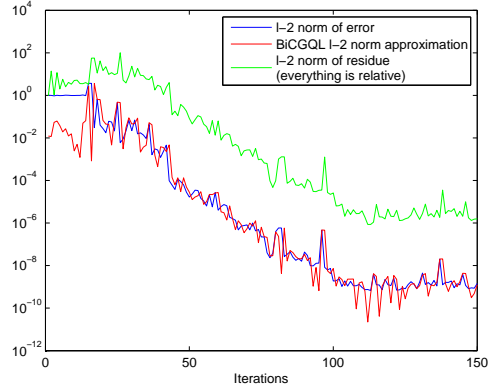}
\par\end{centering}

\caption{Comparison between BiCGQL $l_{2}-$norm estimator, actual $l_{2}-$norm
of the error and $l_{2}-$norm of the residue, condition number of
matrix is $10^{4}$, $d_{1}=d_{2}=0$ \label{fig:Comparison-between-norm}}
\end{figure}

In figure \ref{fig:3}, $A$ is a nonsymmetric matrix. Here, the plot
of approximation vector along with $A-$norm of the error is shown.
Figure \ref{fig:Comparison-between-norm} shows the comparison between
$l_{2}-$norm approximation, actual $l_{2}-$norm of the error and
$l_{2}-$norm of the residue. It is evident that BiCGQL estimators
work efficiently both times. 

For the purpose of extended tests (for both HPD and Indefinite cases),
six different bins of size ten were created with varying condition
number for matrix $A$: 1 to 10, 10 to 100, ... $10^{5}$ to $10^{6}$
etc. For each matrix $A$, 100 different instances of vector $'b'$
were created, each being unique canonical form of order 100. Thus
each bin represents the result accumulated from 1000 different cases.
Below we are comparing our approximation of $A-$norm of error with
estimation by residue vector. Average error in estimating $A-$norm
of error by BiCGQL $A-$norm estimator can be expressed as 

\begin{equation}
\left|\frac{\frac{g_{k}}{||x||_{A}}-\frac{||e_{k}||_{A}}{||x||_{A}}}{\frac{||e_{k}||_{A}}{||x||_{A}}}\right|\label{eq:eq48}
\end{equation}

where $||e||_{A}$ is the $A-$norm of error vector, $g_{k}$ is BiCGQL
$A-$norm estimator, $||r||_{2}$ is the $l_{2}-$norm of residue
vector, $||e||_{2}$ is the $l_{2}-$norm of error vector and $x$
is actual solution vector. While error in estimating $l_{2}-$norm
of the error by residual can be expressed as

\begin{equation}
\left|\frac{\frac{||r_{k}||_{2}}{||b||_{2}}-\frac{||e_{k}||_{2}}{||x||_{2}}}{\frac{||e_{k}||_{2}}{||x||_{2}}}\right|\label{eq:eq49}
\end{equation}

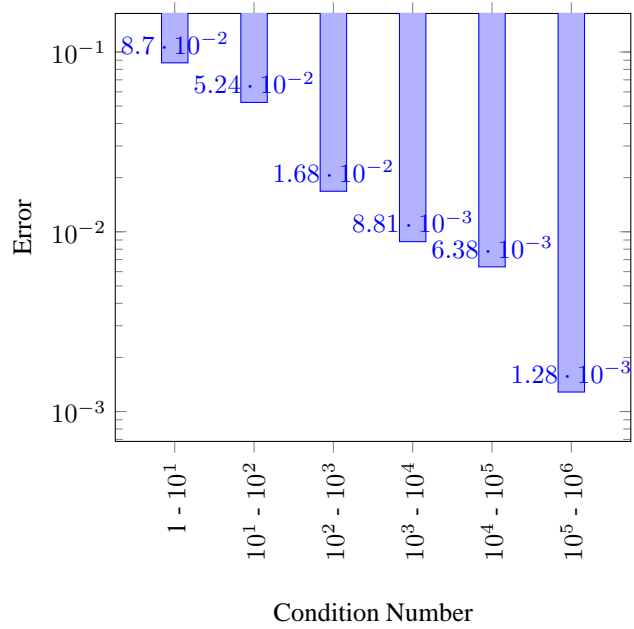
\begin{figure}[h]
\begin{centering}
\begin{tikzpicture} 
\begin{semilogyaxis}[ ybar, enlargelimits=0.15, legend style={at={(0.5,-0.15)}, anchor=north,legend columns=-1},xlabel={Condition Number}, ylabel={Error}, symbolic x coords={1 - $10^1$,$10^1$ - $10^2$,$10^2$ - $10^3$,$10^3$ - $10^4$,$10^4$ - $10^5$,$10^5$ - $10^6$}, xtick=data, nodes near coords, nodes near coords align={vertical}, point meta=rawy , x tick label style={rotate=90,anchor=east},  xlabel style={yshift=-1.5cm},
] 
\addplot coordinates {(1 - $10^1$,0.087) ($10^1$ - $10^2$,0.0524) ($10^2$ - $10^3$,0.0168)($10^3$ - $10^4$,0.00881)($10^4$ - $10^5$,0.00638) ($10^5$ - $10^6$, 0.001284) }; 
\end{semilogyaxis} 
\end{tikzpicture}
\par\end{centering}

\centering{}\caption{Average ratio of relative error in estimating $A-$norm by BiCGQL
and relative error in traditional stopping criteria for an HPD matrix
(each bar shows average over 1000 cases, $d_{1}=d_{2}=4$) \label{fig:HPD-case}}
\end{figure}

\begin{figure}[h]
\begin{centering}
\begin{tikzpicture} 
\begin{semilogyaxis}[ ybar, enlargelimits=0.15, legend style={at={(0.5,-0.15)}, anchor=north,legend columns=-1},xlabel={Condition Number}, ylabel={Error}, symbolic x coords={1 - $10^1$,$10^1$ - $10^2$,$10^2$ - $10^3$,$10^3$ - $10^4$,$10^4$ - $10^5$,$10^5$ - $10^6$}, xtick=data, nodes near coords, nodes near coords align={vertical}, point meta=rawy , x tick label style={rotate=90,anchor=east},  xlabel style={yshift=-1.5cm}
] 
\addplot coordinates {(1 - $10^1$,0.051) ($10^1$ - $10^2$,0.0251) ($10^2$ - $10^3$,0.0082)($10^3$ - $10^4$,0.00594)($10^4$ - $10^5$,0.00218) ($10^5$ - $10^6$, 0.000652) };
\end{semilogyaxis} 
\end{tikzpicture}
\par\end{centering}

\centering{}\caption{Average ratio of relative error in estimating $A-$norm by BiCGQL
and relative error in traditional stopping criteria for a non-HPD
matrix (each bar shows average over 1000 cases, $d_{1}=d_{2}=4$)\label{fig:Indefinite-case}}
\end{figure}
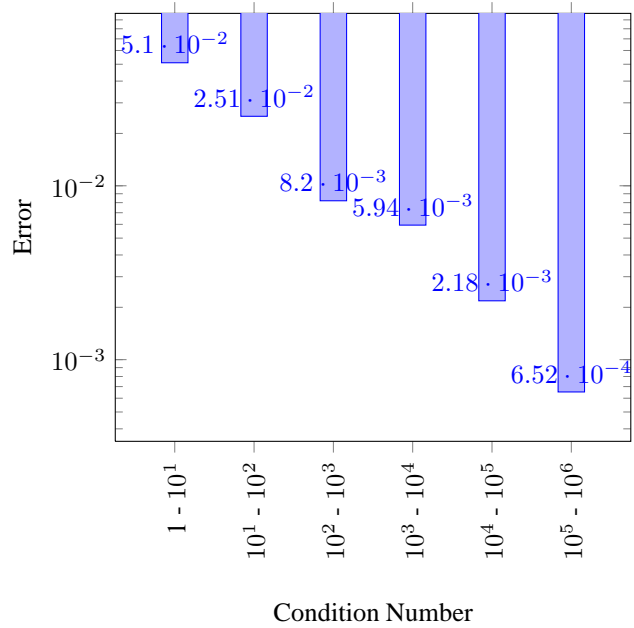

Ratio of equation \ref{eq:eq48} to equation \ref{eq:eq49} would
show the performance of BiCGQL $A-$norm estimator compared to residual
as the estimator of the $l_{2}-$norm of the error. In \ref{fig:HPD-case}
and \ref{fig:Indefinite-case}, each bar represents the average value
of ``ratio of equation \ref{eq:eq48} to equation \ref{eq:eq49}
averaged over all iterations'' over 1000 different cases. Results
obtained for \ref{eq:43} show that the approximation of the $A-$norm
of the error obtained by our $A-$norm estimator is much better than
approximation of the $l_{2}-$norm of the error obtained by the residue
vector. It should be noted that without our approximator iterative
methods would rely on the residue vector which poorly approximated
$l_{2}-$norm of the error. The figures \ref{fig:HPD-case} \& \ref{fig:Indefinite-case}
show that residue keeps becoming unreliable as the condition number
of the problem increases. The graph also shows that our approximator
remains effective in approximating $A-$norm of the error regardless
of the condition number of the problem. 

Similarly, we now show {[}fig. \ref{fig:HPD-case1} \& \ref{fig:Indefinite-case1}{]}
that our $l_{2}-$norm approximator is much better compared to residue
vector in approximating $l_{2}-$norm of the error vector. Test conditions
remain the same as in the previous case. Here, each bar represents
average of following over 1000 different cases... 

\begin{equation}
\left|\frac{\frac{f_{k}}{||x||_{2}}-\frac{||e_{k}||_{2}}{||x||_{2}}}{\frac{||r_{k}||_{2}}{||b||_{2}}-\frac{||e_{k}||_{2}}{||x||_{2}}}\right|
\end{equation}

here $f_{k}$ is BiCGQL $l_{2}-$norm estimator. Even here, BiCGQL
$l_{2}-$norm is proven to be superior than 

\begin{figure}[h]
\begin{centering}
\begin{tikzpicture} 
\begin{semilogyaxis}[ ybar, enlargelimits=0.15, legend style={at={(0.5,-0.15)}, anchor=north,legend columns=-1},xlabel={Condition Number}, ylabel={Error}, symbolic x coords={1 - $10^1$,$10^1$ - $10^2$,$10^2$ - $10^3$,$10^3$ - $10^4$,$10^4$ - $10^5$,$10^5$ - $10^6$}, xtick=data, nodes near coords, nodes near coords align={vertical}, point meta=rawy , x tick label style={rotate=90,anchor=east},  xlabel style={yshift=-1.5cm},
] 
\addplot coordinates {(1 - $10^1$,0.4887) ($10^1$ - $10^2$,0.1839) ($10^2$ - $10^3$,0.0937)($10^3$ - $10^4$,0.00684)($10^4$ - $10^5$,0.00328) ($10^5$ - $10^6$, 0.001434) };
\end{semilogyaxis} 
\end{tikzpicture}
\par\end{centering}

\centering{}\caption{Average ratio of relative error in estimating $l_{2}-$norm by BiCGQL
and relative error in traditional stopping criteria for an HPD matrix
(each bar shows average over 1000 cases, $d_{1}=d_{2}=4$)\label{fig:HPD-case1}}
\end{figure}

\begin{figure}[h]
\begin{centering}
\begin{tikzpicture} 
\begin{semilogyaxis}[ ybar, enlargelimits=0.15, legend style={at={(0.5,-0.15)}, anchor=north,legend columns=-1},xlabel={Condition Number}, ylabel={Error}, symbolic x coords={1 - $10^1$,$10^1$ - $10^2$,$10^2$ - $10^3$,$10^3$ - $10^4$,$10^4$ - $10^5$,$10^5$ - $10^6$}, xtick=data, nodes near coords, nodes near coords align={vertical}, point meta=rawy , x tick label style={rotate=90,anchor=east},  xlabel style={yshift=-1.5cm}
] 
\addplot coordinates {(1 - $10^1$,0.2867) ($10^1$ - $10^2$,0.1578) ($10^2$ - $10^3$,0.0829)($10^3$ - $10^4$,0.00468)($10^4$ - $10^5$,0.00238) ($10^5$ - $10^6$, 0.000434) };
\end{semilogyaxis} 
\end{tikzpicture}
\par\end{centering}

\centering{}\caption{Average ratio of relative error in estimating $l_{2}-$norm by BiCGQL
and relative error in traditional stopping criteria for a non-HPD
matrix (each bar shows average over 1000 cases, $d_{1}=d_{2}=4$)
\label{fig:Indefinite-case1}}
\end{figure}

It is noteworthy that estimator for $l_{2}-$norm of the error holds
greater significance than the estimator for $A-$norm of the error
in most realistic applications. We shall now compare our estimators
with the estimators suggested by Golub and Meurant (\cite{golub2009matrices},
p.210) as below:

\begin{align}
||\epsilon||_{A}^{2} & \approx\frac{(r,Ar)}{(A^{2}r,Ar)}\label{eq:49}\\
||\epsilon||^{2} & \approx\frac{(r,r)^{2}}{(Ar,Ar)}\label{eq:50}
\end{align}

The following expression is averaged over all 1000 different cases,
same as above. (Here $g_{k}^{GM}$ is the estimator suggested by Golub
and Muerant in \ref{eq:49}.)

\begin{equation}
\left|\left(\frac{\frac{g_{k}}{||x||_{A}}-\frac{||e_{k}||_{A}}{||x||_{A}}}{\frac{g_{k}^{GM}}{||x||_{A}}-\frac{||e_{k}||_{A}}{||x||_{A}}}\right)\right|
\end{equation}

\begin{center}
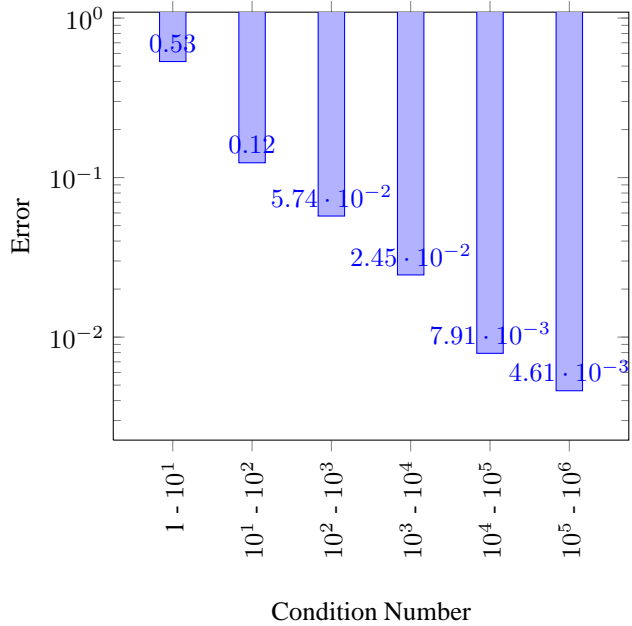
\begin{figure}[H]
\begin{centering}
\begin{tikzpicture} 
\begin{semilogyaxis}[ ybar, enlargelimits=0.15, legend style={at={(0.5,-0.15)}, anchor=north,legend columns=-1},xlabel={Condition Number}, ylabel={Error}, symbolic x coords={1 - $10^1$,$10^1$ - $10^2$,$10^2$ - $10^3$,$10^3$ - $10^4$,$10^4$ - $10^5$,$10^5$ - $10^6$}, xtick=data, nodes near coords, nodes near coords align={vertical}, point meta=rawy , x tick label style={rotate=90,anchor=east},  xlabel style={yshift=-1.5cm}
] 
\addplot coordinates {(1 - $10^1$,0.5337) ($10^1$ - $10^2$,0.1238) ($10^2$ - $10^3$,0.0574)($10^3$ - $10^4$,0.02453)($10^4$ - $10^5$,0.007913) ($10^5$ - $10^6$, 0.004611) };
\end{semilogyaxis} 
\end{tikzpicture}
\par\end{centering}

\centering{}\caption{Average ratio of the relative error in estimating $A-$norm of the
error by BiCGQL and Golub-Meurant estimations for a non-HPD matrix
(each bar shows average over 1000 cases, $d_{1}=d_{2}=0$) \label{fig:compareanorm}}
\end{figure}

\par\end{center}

The following expression is averaged over all 1000 different cases,
same as above. (Here $f_{k}^{GM}$ is the estimator suggested by Golub
and Muerant in \ref{eq:50}.)

\begin{equation}
\left|\left(\frac{\frac{f_{k}}{||x||_{2}}-\frac{||e_{k}||_{2}}{||x||_{2}}}{\frac{f_{k}^{GM}}{||x||_{2}}-\frac{||e_{k}||_{2}}{||x||_{2}}}\right)\right|
\end{equation}

\begin{center}
\begin{figure}[H]
\begin{centering}
\begin{tikzpicture} 
\begin{semilogyaxis}[ ybar, enlargelimits=0.15, legend style={at={(0.5,-0.15)}, anchor=north,legend columns=-1},xlabel={Condition Number}, ylabel={Error}, symbolic x coords={1 - $10^1$,$10^1$ - $10^2$,$10^2$ - $10^3$,$10^3$ - $10^4$,$10^4$ - $10^5$,$10^5$ - $10^6$}, xtick=data, nodes near coords, nodes near coords align={vertical}, point meta=rawy , x tick label style={rotate=90,anchor=east},  xlabel style={yshift=-1.5cm}
] 
\addplot coordinates {(1 - $10^1$,0.4758) ($10^1$ - $10^2$,0.1382) ($10^2$ - $10^3$,0.08223)($10^3$ - $10^4$,0.03371)($10^4$ - $10^5$,0.006771) ($10^5$ - $10^6$, 0.001339) };
\end{semilogyaxis} 
\end{tikzpicture}
\par\end{centering}

\centering{}\caption{Average ratio of the relative error in estimating $l_{2}-$norm of
the error by BiCGQL and Golub-Meurant estimations for a non-HPD matrix(each
bar shows average over 1000 cases, $d_{1}=d_{2}=0$) \label{fig:comparel2norm}}
\end{figure}
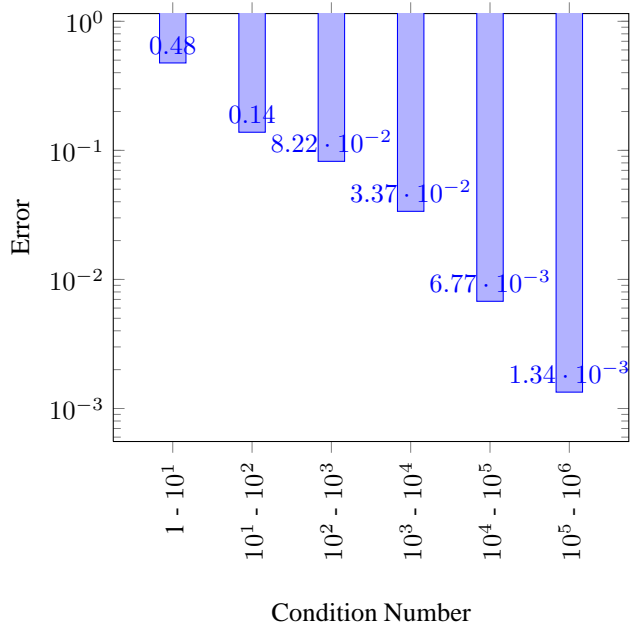

\par\end{center}

From figures \ref{fig:compareanorm} and \ref{fig:comparel2norm},
it can be clearly seen that the estimators in equations \ref{eq:49}
and\ref{eq:50} are good for well-conditioned problem, but they fail
for problems with high conditioned problems. It should also be noted
that these estimators are computationally very expensive as they include
matrix-matrix multiplications. Our estimators (BiCGQL) give much better
results comparatively. Thus it is evident that BiCGQL estimators are
superior in terms of both accuracy and computational cost.

\section{Conclusions}

The importance of BiCGQL estimators are evident for problems with
moderately high condition number $\kappa>100$, and is emphasized
by a few general examples in section 4. The $O(1)$ estimators for
BiCG computations developed by us are on an average $\kappa\times10^{-1}$times
more accurate than residual based stopping criteria and $\kappa\times10^{-2}$
times more accurate than the previously existing estimators. As matrix
$A$ is non-Hermitian, BiCGQL estimators do not necessarily give upper
or lower bounds on the norms of errors, however as previously discussed,
they can be used for indefinite problems. Based on the results presented
in the previous section, we believe that the estimate for the $A-$norm
and $l_{2-}$norm of the error should be implemented into software
realization of the BiCG or similar iterative algorithms as a stopping
criteria instead of the residual. 

\bibliographystyle{IEEEtran}
\bibliography{template}

\end{document}